# ASYMPTOTIC BEHAVIOR OF THE UNCONDITIONAL NPMLE OF THE LENGTH-BIASED SURVIVOR FUNCTION FROM RIGHT CENSORED PREVALENT COHORT DATA[1]

BY MASOUD ASGHARIAN AND DAVID B. WOLFSON

*McGill University*

Right censored survival data collected on a cohort of prevalent cases with constant incidence are length-biased, and may be used to estimate the length-biased (i.e., prevalent-case) survival function. When the incidence rate is constant, so-called stationarity of the incidence, it is more efficient to use this structure for unconditional statistical inference than to carry out an analysis by conditioning on the observed truncation times. It is well known that, due to the informative censoring for prevalent cohort data, the Kaplan–Meier estimator is not the unconditional NPMLE of the length-biased survival function and the asymptotic properties of the NPMLE do not follow from any known result. We present here a detailed derivation of the asymptotic properties of the NPMLE of the length-biased survival function from right censored prevalent cohort survival data with follow-up. In particular, we show that the NPMLE is uniformly strongly consistent, converges weakly to a Gaussian process, and is asymptotically efficient. One important spin-off from these results is that they yield the asymptotic properties of the NPMLE of the incident-case survival function [see Asgharian, M'Lan and Wolfson *J. Amer. Statist. Assoc.* **97** (2002) 201–209], which is often of prime interest in a prevalent cohort study. Our results generalize those given by Vardi and Zhang [*Ann. Statist.* **20** (1992) 1022–1039] under multiplicative censoring, which we show arises as a degenerate case in a prevalent cohort setting.

**1. Introduction.** Left truncated, right censored data have been extensively studied in the statistics literature. (See [4] for a list of references.)

Received October 2001; revised September 2004.
[1]Supported in part by NSERC, Canada.
*AMS 2000 subject classifications.* Primary 62G20; secondary 62N02.
*Key words and phrases.* Prevalent cohort, censored data, length-biased sampling, informative censoring, survival function, NPMLE, consistency, weak convergence and asymptotic efficiency.

Their importance stems from the common use of prevalent cohort study designs to estimate survival from onset of a specified disease (e.g., [13, 25, 26]). In such studies patients are identified with prevalent disease at some instant in calendar time through a cross-sectional survey. Those so identified are then followed forward in time until failure or censoring. Since the possibly censored observed survival times are generated from prevalent cases, they are left truncated. Failure to account for left truncation can lead to substantial overestimation of the survivor function. Indeed, Wolfson et al. [28] showed that survival with dementia from onset had almost certainly been overestimated because of the failure to take left truncation into account. In fact, their adjusted (for left truncation) estimated median survival time was 3.3 years versus 6.6 years for the unadjusted estimated median survival.

When the left truncation time distribution is not specified, the approach to estimating the unbiased survival function is to condition on the observed truncation times (e.g., [1, 27]). However, when there is good reason to assume that the initiation times follow a stationary Poisson process which implies more structure on the truncation times (the so-called *stationarity* assumption, [2]), this special structure may be exploited. Under stationarity, it is not necessary to condition on the observed truncation times and, instead, the natural estimator is the unconditional nonparametric maximum likelihood estimator (NPMLE) [4]. Wang [25] had suggested that an unconditional NPMLE of the unbiased survivor function is more efficient than its conditional counterpart, under stationarity. This improvement in efficiency was later confirmed by Asgharian, M'Lan and Wolfson ([4], Figures 3 and 4, but note the incorrect captions). We shall, in the sequel, reserve the terminology "length-biased" for left truncation under the stationarity assumption.

Several authors have discussed maximum likelihood estimation in the presence of length-biased data [10, 11, 21, 22]. The latter two papers treated the question in the general setting of selection bias, though without allowing for censoring.

This paper establishes the asymptotic properties of the nonparametric maximum likelihood estimator (NPMLE) of the length-biased survival function when the observed data are length-biased and right censored. Now, while the length-biased survival function is itself of little direct interest, the unbiased survival function, which is simply related, is of central importance in a survival analysis based on prevalent cohort data. By exploiting the mapping which relates the length-biased with the unbiased survival function, one may use the results established here to obtain the asymptotic properties of the NPMLE of the unbiased survival function [4]. That is, we present here, for the first time, the foundation upon which the asymptotic inference described in Asgharian, M'Lan and Wolfson rests. It should be pointed out though, that the definition of $\mathcal{G}_2$ given above equation (8) in [4] is erroneous, though with only minor effect on the main result. This error is corrected here.

A subtlety missed by several authors is that the Kaplan–Meier estimator is not suitable as the NPMLE of the length-biased survival function since the right censoring induced by the sampling scheme is informative; in order to be censored, one's failure time must be longer than one's truncation time, that is, be *observable*. (See [4] for further details.)

Our line of attack is similar to that of Vardi [23] and Vardi and Zhang [24], who derived the NPMLE of the length-biased survival function and established its asymptotic properties under multiplicative censoring. They pointed out that the likelihood obtained under multiplicative censoring has the same form as the likelihood obtained from prevalent cohort study data with followup, when the stationarity assumption holds. Importantly, Vardi [23] noted that although the maximum likelihood estimates obtained from these common likelihoods are the same, the asymptotic properties depend on the sampling mechanism that gives rise to the data and must be established afresh in each setting.

It is therefore instructive to place multiplicative censoring in the context of prevalent cohort studies in order to underscore the differences between the two schemes. We also derive the likelihood conditional on the number of censored observations, not because such conditional inference is carried out in practice, but merely to contrast this with the unconditional likelihood. Consider, therefore, the following three situations:

(i) The number of subjects identified at recruitment is $k$. All subjects who are not lost to follow-up are followed until the end of the study period. Those subjects who are lost to follow-up or survive to the end of the study are right censored, the remainder having failed in this time period. It is assumed that the censoring of *the residual life times* (also called *the forward recurrence times*) is random. That is, the times from recruitment until failure are randomly censored. If $M$ denotes the (random) number of uncensored subjects at the end of the study period, then $N = k - M$ denotes the number of censored observations.

(ii) The scenario is the same as that of (i) except that $M$ and $N$ are fixed at the observed values, $m$ and $n$, respectively, and analyses are carried out conditionally.

(iii) The number of subjects identified at the cross-sectional stage is $k = m + n$. At this stage, a fixed number $n$ are immediately censored. The remaining $m$ are followed until failure. It is easily seen that this sampling scheme is equivalent to that of multiplicative censoring.

The setup described by (i) occurs in practice most frequently, and is the focus of this paper. In fact, in Section 7 we show explicitly how sampling scheme (iii) arises as a particular, although degenerate, case of scheme (i).

Under (iii) censoring is precluded after recruitment, which is clearly an unrealistic assumption in a prevalent cohort study with follow-up. See [23] though for examples of multiplicative censoring in different contexts. Section 7 contains further discussion on schemes (ii) and (iii).

The above generalities are perhaps better understood through The Canadian Study of Health and Aging (CSHA), a large prevalent cohort study with follow-up conducted to investigate, primarily, various aspects of dementia in the elderly Canadian population.

Briefly, during a six-month period in 1991 roughly 10,000 Canadians over the age of 65 were recruited and screened for prevalent dementia. Dementias considered included mainly probable Alzheimer's disease, possible Alzheimer's disease and vascular dementia. At the time of diagnosis, age at onset was ascertained from the patient's caregiver. Those subjects diagnosed with dementia in 1991 were followed until censoring or death. Follow-up ended in 1996, and subjects who were still alive were deemed to have been right censored. Very few subjects were lost to follow-up (also considered to be right censored) between 1991 and 1996. Times of death from any cause, or of censoring, were recorded for all subjects diagnosed with dementia. (See [28] for further details.)

One of the many aims of the CSHA was to estimate the unbiased survival function of subjects with incident dementia, where the origin was date of onset and the endpoint was death from any cause. The data available for this estimation problem had several features: (a) They were left truncated because subjects with dementia were identified as prevalent rather than incident cases. (b) The underlying process that generated the onset times of dementia was thought to have been roughly stationary in the sense that the intensity function of the initiation process was constant; the basis for this assumption is discussed by Asgharian, M'Lan and Wolfson ([4], Figure 5 (note the incorrect caption)) and Asgharian, Wolfson and Zhang [2]. The observed survival times were, therefore, length-biased. (c) The interval from onset to recruitment, *the current life time* (also called *the backward recurrence time*), as well as the minimum of the interval from recruitment to death, the residual life time and censoring, were recorded. That is, there was more information than that contained in the length-biased, possibly censored, survival times alone. (d) The censoring of the residual life times was random.

Inference about the length-biased survival function from onset with dementia could be based on the current life times together with their possibly randomly censored residual life times, as was done by Asgharian, M'Lan and Wolfson [4]. Alternatively, with the same data, though less conventionally, one could have conditioned on the observed number of censored subjects and those who died, so that all ensuing inference would have been conditional [see (ii)].

It is of interest to note that length-biased sampling also arises when a renewal process is sampled at some point $t$. The interval surrounding $t$ has the density $g(x) = xf(x)/\mu_X$, where $f$ is the density and $\mu_X$ the (finite) mean of the sojourn time distribution. See [2] for the differences between the two settings and a characterization of stationarity. Other sampling schemes have been discussed in the literature and can be depicted using the Lexis diagram [7, 17, 18], and it is possible to carry out inference from a prevalent cohort study under appropriate restrictions when there is no follow-up [15, 16, 20].

The layout of the paper is as follows: In Section 2 we present the likelihoods for the different sampling schemes. Thereafter, we focus on sampling scheme (i). A general overview of the proofs is given in Section 3. In Section 4 we establish uniform consistency of the NPMLE, and in Section 5 we discuss weak convergence of the NPMLE. Asymptotic efficiency of the NPMLE is presented in Section 6. In Section 7 we return briefly to schemes (ii) and (iii), expanding on their relationships to scheme (i). Section 8 summarizes our results and contains some concluding comments.

## 2. Preliminaries and the likelihoods for different sampling schemes.

2.1. *Preliminaries.* While "stationarity" refers to the pattern of chronological cross-sectional sampling, all of the notation below relates to current-age (current life time), failure and censoring durations for the individual subjects sampled. Suppose that associated with each subject in a target population we have a triple $(X', T', C')$, where $X'$ represents the failure time, $T'$ the truncation time and $C'$ the censoring time. Often a reasonable assumption is that $X'$ is independent of $(T', C')$, while $P(C' \geq T') = 1$ [25]. In a cross-sectional survey subjects are observed only if $X' \geq T'$. Under the stationarity assumption, the survival time density of the *observed* subjects, the length-biased density, say $g$, is related to $f_{X'}$, the unbiased density, through the equation

$$g(x) = f_{X'|X' \geq T'}(x|X' \geq T') = \frac{xf_{X'}(x)}{\mu_{X'}}.$$

We give, initially, the likelihoods derived under (i), (ii) and (iii), in order to emphasize the differences between the three situations. Although the likelihoods are similar, the appearance of the random censoring indicator under scheme (i) ($M$ random) requires special treatment in the derivation of the large sample properties of the maximum likelihood estimator. These properties are established in Sections 4, 5 and 6 under scheme (i).

None of the three likelihoods below depends explicitly on the residual and current life times separately. However, the derivations of these likelihoods depend explicitly on knowledge of the current life times, as well as the randomly right censored residual life times. These are the data typically observed in a prevalent cohort study with follow-up, when the times of onset are known.

Associated with each *observed* subject in a prevalent cohort study, we have a triple,

$$(A_i, R_i \wedge C_i, \delta_i), \qquad i = 1, 2, \ldots, k,$$

where $A_i$, $R_i$ and $C_i$ are, respectively, the current-age, the residual life time and the residual censoring time for the $i$th *observed* subject. The indicator function $\delta_i$ is the censoring indicator, that is,

$$\delta_i = \begin{cases} 1, & \text{if the } i\text{th subject is not censored } (R_i \leq C_i), \\ 0, & \text{otherwise.} \end{cases}$$

It is reasonable in many applications to assume that $C_i$ is independent of $(A_i, R_i)$. We adopt this assumption in the sequel. The vectors $(A_i, R_i \wedge C_i, \delta_i)$, $i = 1, 2, \ldots, k$, are also assumed to be independent.

Note that the failure time and censoring time associated with each *observed* subject are, respectively, $X' = A + R$ and $Y' = A + C$. One can therefore easily show that, in general, if $C$ is independent of $(A, R)$, then $\text{Cov}(X, Y) = \sigma_A^2 [1 + \rho_{A,R} \sigma_R / \sigma_A]$, where $\sigma_A^2 = \text{Var}(A)$, $\sigma_R^2 = \text{Var}(R)$ and $\rho_{A,R} = \text{corr}(A, R)$. Thus, except for trivial cases, failure times and censoring times are positively correlated under stationarity, since stationarity implies $A$ is conditionally $\text{Unif}(0, X')$ given $X'$, so that $\sigma_A = \sigma_{X'-A} = \sigma_R$. This then implies that the censoring mechanism in the setting under study is informative.

Under the stationarity assumption,

$$(2.1) \qquad f_{A,R}(a, r) = \begin{cases} \dfrac{f_{X'}(a+r)}{\mu_{X'}}, & \text{if } a, r > 0, \\ 0, & \text{otherwise,} \end{cases}$$

which corresponds to the well-known expression for the joint density of the current and residual life times, respectively, of a renewal process (see [9, 23]). In the sequel we use $f_U$ for $f_{X'}$, where $U$ in the subscript stands for "unbiased." Using (2.1), one can easily derive the distribution function, say $G$, of $X = A + R$, the length-biased survival time, whose density is

$$(2.2) \qquad g(x) = \frac{x f_U(x)}{\mu_U}.$$

Let

$$G_*(t) = P(A + R \leq t | \delta = 1),$$

with density function $g_*(t)$. We then have

$$g_*(t) = \frac{1}{p} \int_0^t f_{A,R}(t-r, r) S_C(r) \, dr$$

$$\begin{aligned}
&= \frac{1}{p} \int_0^t \frac{f_U(t)}{\mu_U} S_C(r) \, dr \\
&= \frac{f_U(t)}{p\mu_U} \int_0^t S_C(r) \, dr \\
&= \frac{g(t)}{pt} \int_0^t S_C(r) \, dr,
\end{aligned}$$
(2.3)

where $p = P(\delta = 1) = P(R \leq C)$ and $S_C(r) = 1 - F_C(r) = 1 - P(C \leq r)$.

Suppose $F_*(t) = P(A + C \leq t | \delta = 0)$, with density function $f_*(t)$. Then

$$\begin{aligned}
f_*(t) &= \frac{1}{1-p} \int_0^t \int_c^\infty f_{A,R}(t - c, r) \, dr \, dF_C(c) \\
&= \frac{1}{1-p} \int_0^t \int_c^\infty \frac{f_U(t + r - c)}{\mu_U} \, dr \, dF_C(c) \\
&= \frac{1}{(1-p)\mu_U} \int_0^t \int_t^\infty f_U(u) \, du \, dF_C(c) \\
&= \frac{S_U(t) F_C(t)}{\mu_U (1-p)} = \frac{f(t) F_C(t)}{1-p},
\end{aligned}$$
(2.4)

where

$$f(t) = \frac{S_U(t)}{\mu_U} = \int_t^\infty z^{-1} \, dG(z)$$

is the residual lifetime density. We turn now to the derivation of the likelihoods under the schemes (i), (ii) and (iii) of Section 1.

2.2. *Random censoring ($M$ random).* This is the case in which $M$ and $N = k - M$ are random and arises under situation (i) of Section 1. The observations comprise

$$(A_i, R_i \wedge C_i, \delta_i), \qquad i = 1, 2, \ldots, k.$$

Let $\mathcal{UC}$ and $\mathcal{C}$ denote, respectively, the sets of indices of the uncensored and censored observations. Let $a_i$, $r_i$ and $c_i$ denote, respectively, the realized values of $A_i$, $R_i$ and $C_i$, and let $x_i = a_i + r_i$, $y_j = a_j + c_j$, and $z = a_j + r$, for $i \in \mathcal{UC}$ and $j \in \mathcal{C}$. The likelihood is

$$\begin{aligned}
\mathcal{L}_\mathcal{R} &= \left( \prod_{i \in \mathcal{UC}} f_{A,R}(a_i, r_i) \right) \left( \prod_{j \in \mathcal{C}} dP(a_j, R_j \geq c_j) \right) \\
&\stackrel{\text{using (2.1)}}{\propto} \left( \prod_{i \in \mathcal{UC}} dG(x_i) \right) \left( \prod_{j \in \mathcal{C}} \int_{c_j \leq r} \frac{f_U(a_j + r)}{\mu_U} \, dr \right)
\end{aligned}$$
(2.5)

$$= \left(\prod_{i \in \mathcal{UC}} dG(x_i)\right)\left(\prod_{j \in \mathcal{C}} \int_{y_j \leq z} z^{-1}\, dG(z)\right)$$

$$= \prod_{i=1}^{k} (dG(x_i))^{\delta_i} \left(\int_{y_i \leq z} z^{-1}\, dG(z)\right)^{1-\delta_i},$$

which has a form different from the likelihood that leads to the Kaplan–Meier estimator in the presence of randomly right censored data. (See, e.g., [14], page 15.)

2.3. *Random censoring (conditional on M)*. Here the data arise as in Section 2.2, but all analyses are carried out conditional on $M = m$ and $N = k - m = n$. The "effective" observations comprise, therefore,

$$(A_i, R_i) \sim f_{A,R|\delta=1}, \qquad i = 1, 2, \ldots, m,$$

and

$$(A_j, C_j) \sim f_{A,C|\delta=0}, \qquad j = 1, 2, \ldots, n.$$

The likelihood contributions are

$$f_{A,R}(A=a, R=r|\delta=1) = f_{A,R}(a, r|R \leq C)$$
$$= \frac{1}{p} S_C(r) f_{A,R}(a, r)$$
$$= \frac{S_C(r)}{p} \frac{f_U(a+r)}{\mu_U}$$
$$= \frac{S_C(r)}{p(a+r)} g(a+r)$$

and

$$f_{A,C}(A=a, C=c|\delta=0) = f_{A,C}(a, c|R \geq C)$$
$$= \frac{1}{1-p} f_C(c) \int_c^\infty f_{A,R}(a, r)\, dr$$
$$= \frac{f_C(c)}{1-p} \int_c^\infty \frac{f_U(a+r)}{\mu_U}\, dr$$
$$= \frac{f_C(c)}{1-p} \int_{a+c \leq z} z^{-1}\, dG(z).$$

Thus, the likelihood

$$\mathcal{L}_\mathcal{S} = \left(\prod_{i=1}^{m} dG(x_i) \frac{S_C(r_i)}{p x_i}\right)\left(\prod_{j=1}^{n} \int_{y_j \leq z} z^{-1}\, dG(z) \frac{f_C(c_j)}{1-p}\right)$$

$$\propto \left(\prod_{i=1}^{m} dG(x_i)\right)\left(\prod_{j=1}^{n} \int_{y_j \leq z} z^{-1} dG(z)\right).$$

2.4. *Multiplicative censoring.* The scenario described in (iii) of Section 1 is equivalent to Vardi's [23] scheme of multiplicative censoring. In the context of a prevalent cohort study, multiplicative censoring is induced by defining the distribution function of $C$, the residual censoring time, as

$$(2.6) \qquad F_C(t) = \begin{cases} 0, & \text{if } t < 0, \\ 1-p, & \text{if } 0 \leq t < \tau, \\ 1, & \text{if } t \geq \tau, \end{cases}$$

where $\tau = \inf\{t : G(t) = 1\}$. If the residual censoring distribution is given by (2.6), then $f_* \equiv f$ and $g_* \equiv g$. Vardi and Zhang [24] considered a sequence of $\{F_C^k\}$ with $p_k = m/k$, so that they essentially conditioned on the censoring proportion.

**3. Asymptotics: general overview and master equation.** The discussion in Sections 3, 4, 5 and 6 is restricted to sampling scheme (i). Let $\widehat{G}$ be the maximizer of the likelihood $\mathcal{L}_\mathcal{R}$ [see (2.5)] with respect to $G$. In this section we present the master equation, and outline the main steps in establishing the uniform strong consistency of $\widehat{G}$ and show that $U_{m,n}$, defined below, converges weakly to a Gaussian process. The details are given in Sections 4, 5 and the Appendix.

Define
$$U_{m,n} = \sqrt{k}(\widehat{G} - G),$$
$$W_{X,m} = \sqrt{m}(G_m - G_*)$$

and

$$W_{Y,n} = \sqrt{n}(F_n - F_*),$$

where $m = \sum_{i=1}^{k} \delta_i$ and $n = k - m$ are realized values of $M$ and $N$, respectively, and $G_m$ and $F_n$ are the respective empirical distribution functions of $x_1, \ldots, x_m$ and $y_1, \ldots, y_n$. Let $\hat{p} = \hat{p}_k = \frac{m}{k}$ and let $t_1 < \cdots < t_h$ be the distinct values of $x_1, \ldots, x_m$ and $y_1, \ldots, y_n$.

The derivation of the asymptotics begins with the score equation derived from the likelihood $\mathcal{L}_\mathcal{R}$. The NPMLE must satisfy the score equation

$$(3.1) \qquad d\widehat{G}(t) = \hat{p}\, dG_m(t) + (1-\hat{p}) \int_{0 < y \leq t} \frac{dF_n(y)}{\int_{y \leq z} z^{-1} d\widehat{G}(z)} t^{-1} d\widehat{G}(t),$$

subject to $\sum_{j=1}^{h} d\widehat{G}(t_j) = 1$ and $d\widehat{G}(t_j) > 0$, $j = 1, \ldots, h$ ([23], page 754). Integrating both sides of (3.1), we obtain

$$(3.2) \quad \widehat{G}(t) = \hat{p} G_m(t) + (1-\hat{p}) \int_{0 < x \leq t} \left[\int_{0 < y \leq x} \frac{dF_n(y)}{\int_{y \leq z} z^{-1} d\widehat{G}(z)}\right] x^{-1} d\widehat{G}(x),$$

where the final integrand is defined to be 0 for $x > t_h$.

Our first objective is to use (3.2) to provide an explicit linear mapping on a function space $D_0[0, t]$ (see Section 4 for precise definitions) expressing an explicit linear functional of $U_{m,n}$ approximately as a linear functional of $W_{X,m}$, $W_{Y,n}$ and $\hat{p} - p$. The linear functional of $U_{m,n}$ is shown to be boundedly invertible, and the resulting expression for $U_{m,n}$ is used to prove uniform consistency and efficiency for $\widehat{G}$ and weak distributional convergence for $U_{m,n}$.

LEMMA 1 (Master equation). *Let*

$$\hat{f}(t) = \int_{t<z} z^{-1} \, d\widehat{G}(z),$$

(3.3)
$$W_{m,n}(t) = \hat{p}^{1/2} W_{X,m}(t) + (1-\hat{p})^{1/2} \hat{f}(t) \int_{0<y\leq t} W_{Y,n}(y) \, d\left(\frac{1}{\hat{f}(y)}\right)$$

*and*

(3.4) $\quad V_{m,n}(t) = W_{m,n}(t) + p\left(\frac{p}{1-p}\right)^{1/2} (G_*(t) - G(t)) \dfrac{\sqrt{k}(\hat{p} - p)}{\sqrt{p(1-p)}}.$

*Then*

$$\left(\frac{\hat{p} - p}{1 - p}\right) U_{m,n}(t)$$

$$+ \left[1 - \left(\frac{\hat{p} - p}{1 - p}\right)\right]$$

(3.5) $\quad\quad \times \left\{ p \int_{0<x\leq t} \frac{g_*(x)}{g(x)} \, dU_{m,n}(x) \right.$

$$\left. + (1-p) \int_{0<y\leq t} y \left( \int_{y\leq z} \frac{U_{m,n}(z)}{z^2} \, dz \right) d\left[\left(\frac{\hat{f}(t)}{\hat{f}(y)} - 1\right)\frac{f_*(y)}{f(y)}\right] \right\}$$

$$= V_{m,n}(t),$$

*where* $m = \sum_{i=1}^{k} \delta_i$, $n = k - m$ *and* $p = P(\delta = 1) = P(R \leq C)$.

PROOF. See the Appendix. □

Lemma 1 relates $U_{m,n}$ to the empirical processes $W_{X,m}$ and $W_{Y,n}$, which are indexed by the realized values of the random integers $M$ and $N$. Equation (A.5) shows that the process $V_{m,n}(t)$, given by (3.4), can be expressed as the image of a linear operator applied to $U_{m,n}$. To see this, define

(3.6) $\quad \mathcal{G}_{k,1}(u)(t) = p\left(\frac{1-\hat{p}}{1-p}\right) \int_{0<x\leq t} \frac{g_*(x)}{g(x)} \, du(x)$

$$\mathcal{G}_{k,2}(u)(t) = (1-p)\left(\frac{1-\hat{p}}{1-p}\right)$$
(3.7)
$$\times \int_{0<y\leq t} y\left(\int_{y\leq z} \frac{u(z)}{z^2}\,dz\right) d\left[\left(\frac{\hat{f}(t)}{\hat{f}(y)}-1\right)\frac{f_*(y)}{f(y)}\right]$$

and

(3.8) $$\mathcal{H}_k(u)(t) = \left(\frac{\hat{p}-p}{1-p}\right)u(t).$$

Define

$$\mathcal{G}_k = \mathcal{G}_{k,1} + \mathcal{G}_{k,2},$$

and express $\mathcal{F}_k$ as

$$\mathcal{F}_k = \mathcal{H}_k + \mathcal{G}_k.$$

Then we may write

(3.9) $$\mathcal{F}_k(U_{m,n}) = V_{m,n}.$$

It is clear that $\mathcal{G}_{k,1}$, $\mathcal{G}_{k,2}$, $\mathcal{H}_k$ and, thus, $\mathcal{F}_k$, are linear operators.

**4. Uniform consistency of $\widehat{G}$.** To study the properties of $\mathcal{F}_k$, we first need to determine the space on which $\mathcal{F}_k$ acts. Let $D_0[0,t]$ be the space of all cadlag functions $u(\cdot)$ on $[0,t]$ that vanish at 0. The space $D_0[0,t]$ endowed with the uniform topology, the topology induced by the uniform norm, $\|u\| = \sup_{0\leq s<t}|u(s)|$, is a Banach space. This implies that $\mathcal{L}(D_0[0,t], D_0[0,t])$, the space of bounded linear operators on $D_0[0,t]$, is a Banach algebra. The other fact about $D_0[0,t]$, endowed with the uniform topology, that we need in the sequel is that cadlag functions have countably many jumps [19]. This guarantees that cadlag functions are Riemann integrable on bounded intervals.

Define $\tau = \inf\{t : G(t) = 1\}$. Let $\tau < \infty$ and

$$\mathcal{J} = \left\{t \leq \tau : \left(\frac{2}{\alpha(t)} - \frac{1}{1-\beta}\right)\beta < 1\right\},$$

where

$$\alpha(t) = \frac{1}{t}\int_0^t S_C(s)\,ds = \frac{pg_*(t)}{g(t)} \quad \text{and} \quad \beta = F_C(0) > 0.$$

We note that $\alpha(t)$ is a decreasing function with $\lim_{t\to 0}\alpha(t) = 1 - \beta$. Thus, $\beta < 1/2$ if $t \in \mathcal{J}$. It is also easy to see that a sufficient condition for $t \in \mathcal{J}$ is $F_C(t) < 1/2$, since $\alpha(t) \geq S_C(t)$. For the interpretation of this condition, see [4]. The condition on $t$ is somewhat less restrictive than that given in [4]. See also Section 8 for further comments.

LEMMA 2. *If $t \in \mathcal{J}$, then a.s. for all sufficiently large $k$,*

(a) $\mathcal{F}_k$ *is a bounded linear operator on $D_0[0,t]$, and*

$$\|\mathcal{F}_k\| \leq \frac{|\hat{p}-p|}{1-p} + \frac{1-\hat{p}}{1-p}; \tag{4.1}$$

(b) $\mathcal{F}_k$ *is an invertible linear operator on $D_0[0,t]$, and*

$$\|\mathcal{F}_k^{-1}\| \leq \frac{\widehat{\lambda}(t)}{1-\widehat{\lambda}(t)|\hat{p}-p|/(1-p)}, \tag{4.2}$$

*where*

$$\widehat{\lambda}(t) = \frac{((1-p)/(1-\hat{p}))(2/\alpha(t) - 1/(1-\beta))}{1 - (2/\alpha(t) - 1/(1-\beta))\beta}.$$

PROOF. (a) Define $A_k(u)(s) = \hat{f}(s) \int_0^s y \int_{y \leq z} \frac{u(z)}{z^2} dz\, d\frac{1}{\hat{f}(y)}$. Then $\|A_k\| \leq 1$ and, therefore, $\|\mathcal{G}_{k,2}\| \leq \beta \frac{1-\hat{p}}{1-p}$, via (2.4). On the other hand, using integration by parts and since

$$\alpha(x) = \frac{pg_*(x)}{g(x)} \to \begin{cases} 1-\beta, & \text{as } x \to 0, \\ 0, & \text{as } x \to \infty, \end{cases}$$

we have $\|\mathcal{G}_{k,1}\| \leq (1-\beta)\frac{1-\hat{p}}{1-p}$. This completes the proof of part (a).

(b) We have that $\mathcal{G}_{k,1}$ is invertible and

$$\mathcal{G}_{k,1}^{-1}(u)(s) = \frac{1-p}{1-\hat{p}} \int_{0 \leq x \leq s} \frac{g(x)}{pg_*(x)}\, du(x).$$

Using integration by parts and as $\alpha(x)$ is a decreasing function, we have $\|\mathcal{G}_{k,1}^{-1}\| \leq \frac{1-p}{1-\hat{p}}(\frac{2}{\alpha(t)} - \frac{1}{1-\beta})$ a.s. Since $\mathcal{L}(D_0[0,t], D_0[0,t])$ is a Banach algebra, $\mathcal{G}_k$ is invertible a.s. for large $k$. In fact,

$$\mathcal{G}_k = \mathcal{G}_{k,1}(I + \mathcal{G}_{k,1}^{-1}\mathcal{G}_{k,2})$$

and thus

$$\mathcal{G}_k^{-1} = (I + \mathcal{G}_{k,1}^{-1}\mathcal{G}_{k,2})^{-1}\mathcal{G}_{k,1}^{-1},$$

which implies that

$$\|\mathcal{G}_k^{-1}\| \leq \|(I + \mathcal{G}_{k,1}^{-1}\mathcal{G}_{k,2})^{-1}\|\|\mathcal{G}_{k,1}^{-1}\|$$
$$\leq \frac{((1-p)/(1-\hat{p}))(2/\alpha(t) - 1/(1-\beta))}{1 - (2/\alpha(t) - 1/(1-\beta))\beta} = \widehat{\lambda}(t) \quad \text{a.s.,}$$

since $\hat{p} \to p$ a.s. as $k$ gets large. Having established the invertibility of $\mathcal{G}_k$, we have

$$\mathcal{F}_k = \mathcal{G}_k(I + \mathcal{G}_k^{-1}\mathcal{H}_k).$$

Using the strong consistency of $\hat{p}$ and the fact that $\mathcal{L}(D_0[0,t], D_0[0,t])$ is a Banach algebra, we obtain once again

$$\mathcal{F}_k^{-1} = (I + \mathcal{G}_k^{-1}\mathcal{H}_k)^{-1}\mathcal{G}_k^{-1} \quad \text{a.s.}$$

and also

$$\|\mathcal{F}_k^{-1}\| \leq \|(I + \mathcal{G}_k^{-1}\mathcal{H}_k)^{-1}\|\|\mathcal{G}_k^{-1}\|$$

$$\leq \frac{\|\mathcal{G}_k^{-1}\|}{1 - \|\mathcal{G}_k^{-1}\||\hat{p}-p|/(1-p)} \leq \frac{\widehat{\lambda}(t)}{1 - \hat{\lambda}(t)|\hat{p}-p|/(1-p)} \quad \text{a.s.}$$

This completes the proof for part (b). □

Theorem 1 below and its corollary prove the uniform strong consistency of $\widehat{G}$.

THEOREM 1. *Let $\widehat{G}$ be the NPMLE of the continuous lifetime distribution function $G$, and $t \in \mathcal{J}$. Then*

$$(4.3) \quad \|\widehat{G} - G\|_\infty = \sup_{0 \leq s \leq t} |\widehat{G}(s) - G(s)| = O\left(\sqrt{\frac{\log \log k}{k}}\right) \quad a.s.$$

PROOF. Using Lemma 2 and (3.9),

$$\|\widehat{G} - G\|_\infty \leq \|\mathcal{F}_k^{-1}\|\left\|\frac{V_{m,n}}{\sqrt{k}}\right\|.$$

On the other hand,

$$\limsup_{k \to \infty} \|\mathcal{F}_k^{-1}\| \leq \lambda(t) = \frac{2/\alpha(t) - 1/(1-\beta)}{1 - (2/\alpha(t) - 1/(1-\beta))\beta} \quad \text{as } k \to \infty.$$

It therefore suffices to show that

$$(4.4) \quad \left\|\frac{V_{m,n}}{\sqrt{k}}\right\|_\infty = O\left(\sqrt{\frac{\log \log k}{k}}\right) \quad \text{a.s.}$$

Next, using (3.3), we have

$$|W_{m,n}(t)| \leq \hat{p}^{1/2}|W_{X,m}(t)| + (1-\hat{p})^{1/2}\|W_{Y,n}\|_\infty.$$

On the other hand, using the law of the iterated logarithm (LIL),

$$\left\|p\left(\sqrt{\frac{p}{1-p}}\right)(G_*(t) - G(t))\frac{\hat{p}-p}{\sqrt{p(1-p)}}\right\|_\infty = O\left(\sqrt{\frac{\log \log k}{k}}\right) \quad \text{a.s.}$$

To complete the proof, we need to show that

$$(4.5) \quad \|G_m - G_*\|_\infty = O\left(\sqrt{\frac{\log \log k}{k}}\right) \quad \text{a.s.}$$

and

$$\|F_n - F_*\|_\infty = O\left(\sqrt{\frac{\log \log k}{k}}\right) \quad \text{a.s.} \tag{4.6}$$

To establish (4.5), one may either use the Kolmogorov exponential bounds and Borel–Cantelli lemma or, by using the LIL, argue as follows. Splitting one of the sums into two parts, we have that

$$\|G_m - G_{kp}\|_\infty \leq \frac{1}{kp}\left[|[kp] - m| + mO\left(\sqrt{\frac{\log \log k}{k}}\right)\right],$$

which implies

$$\|G_m - G_{kp}\|_\infty = O\left(\sqrt{\frac{\log \log k}{k}}\right) \quad \text{a.s.,} \tag{4.7}$$

where $G_{kp}$ is the empirical distribution function of $x_1, x_2, \ldots, x_{[kp]}$. Now, using the triangle inequality,

$$\|G_m - G_*\|_\infty \leq \|G_m - G_{kp}\|_\infty + \|G_{kp} - G_*\|_\infty. \tag{4.8}$$

Thus, (4.5) follows from (4.7), (4.8) and the fact that

$$\|G_{kp} - G_*\|_\infty = O\left(\sqrt{\frac{\log \log k}{k}}\right) \quad \text{a.s.}$$

Likewise, one can establish (4.6). This completes the proof. □

Equation (4.3), established by Theorem 1, tells us how fast $\widehat{G}$ converges to $G$ in the supnorm topology. Strong consistency of $\widehat{G}$ may therefore be stated as a corollary to Theorem 1.

COROLLARY 1. *Suppose $F_U$ is a continuous lifetime distribution function. Let $G$ be the length-biased distribution function of $F_U$ given by* (2.2) *and $t \in \mathcal{J}$. Then $\widehat{G}$, the NPMLE of $G$ based on data collected according to sampling scheme* (i), *is uniformly strongly consistent on $[0, t]$.*

**5. Weak convergence of $U_{m,n}$.** To establish the weak convergence of $U_{m,n}$ to a Gaussian process, we first need to prove the following lemma.

LEMMA 3. *If $t \in \mathcal{J}$, then*

$$\|\mathcal{F}_k(u) - \mathcal{F}(u)\|_\infty \to 0 \quad a.s. \ \forall \, u \in D_0[0, t],$$

*where*

(5.1) $$\mathcal{F} = \mathcal{G}_1 + \mathcal{G}_2,$$

(5.2) $$\mathcal{G}_1(u)(s) = p \int_{0 < x \leq s} \frac{g_*(x)}{g(x)} \, du(x)$$

*and*

(5.3) $$\mathcal{G}_2(u)(s) = (1-p) \int_{0 < y \leq s} y \left( \int_{y \leq z} \frac{u(z)}{z^2} \, dz \right) d\left[ \left( \frac{f(s)}{f(y)} - 1 \right) \frac{f_*(y)}{f(y)} \right].$$

PROOF. Using the law of large numbers and the bound
$$\|\mathcal{H}_k\| \leq \frac{|\hat{p} - p|}{1 - p} \quad \text{a.s.},$$
we have $\|\mathcal{H}_k\| \to 0$ a.s. as $k \to \infty$. It is also easily seen that
$$\|\mathcal{G}_{k,1} - \mathcal{G}_1\| \leq \left| \frac{1 - \hat{p}}{1 - p} - 1 \right| \to 0 \quad \text{a.s. as } k \to \infty.$$
To complete the proof, we need to show that
$$\|\mathcal{G}_{k,2}(u) - \mathcal{G}_2(u)\|_\infty \to 0 \quad \text{a.s. as } k \to \infty$$
for all $u \in D_0[0,t]$. This can be done along the lines of Lemma 2 of [24]. We therefore omit the proof. □

THEOREM 2. *Suppose $F_U$ is a continuous life time distribution function. Let $G$ be the length-biased distribution function of $F_U$ given by* (2.2) *and let $\widehat{G}$ be the NPMLE of $G$. Then for any $t \in \mathcal{J}$,*
$$U_{m,n} = \sqrt{k}(\widehat{G} - G) \xrightarrow{w} U = \mathcal{F}^{-1}(V) \quad \text{in } D_0[0,t],$$
*where $\mathcal{F}^{-1}$ is the inverse of $\mathcal{F}$ given by* (5.1), (5.2) *and* (5.3),
$$V(s) = p^{1/2} B_1(G_*(s)) + (1-p)^{1/2} f(s) \int_{0 < y \leq s} B_2(F_*(y)) \, d\frac{1}{f(y)}$$
$$+ p \left( \frac{p}{1-p} \right)^{1/2} (G_*(s) - G(s)) Z,$$
*$Z \sim N(0,1)$, and where $B_1$ and $B_2$ are independent Brownian bridge processes, independent of $Z$.*

PROOF. We first need to show that $\mathcal{F}$ is invertible. This is done using a similar argument to that used in Lemma 2. It is also easy to see that $\|\mathcal{F}^{-1}\| \leq \lambda(t)$, where $\lambda(t) = \frac{(2/\alpha(t) - 1/(1-\beta))}{1 - (2/\alpha(t) - 1/(1-\beta))\beta}$. In view of Theorem 7.3.2 and 7.3.3 of [8], the limiting processes of $W_{m,n}$ and $W_{[kp],[k(1-p)]}$ are the same. Hence, we need only find the limiting process of $W_{[kp],[k(1-p)]}(s) + p(\frac{p}{1-p})^{1/2}(G_*(s) - G(s))\frac{\sqrt{k}(\hat{p}-p)}{\sqrt{p(1-p)}}$. Theorem 2 of [24] may now be used to complete the proof. $\square$

**6. Asymptotic efficiency of $\widehat{G}$.** It transpires that, under scheme (i), $\widehat{G}$ is asymptotically efficient in the class of regular estimators whose finite-dimensional limiting laws are continuous in $G$. This result is perhaps not unexpected, given that Vardi and Zhang [24] have established the asymptotic efficiency of $\widehat{G}$ under multiplicative censoring, which we show in Section 7 is a special case of scheme (i). Since the proof under scheme (i) mimics that of Vardi and Zhang, we omit the details. A systematic account of asymptotic efficiency and the convolution theorem can be found in [6, 12]. Here we follow the approach taken by Vardi and Zhang [24] and confine our attention to regular estimators whose finite-dimensional limiting laws are continuous in $G$.

Let $H(\cdot)$ be a stochastic process in $D_0[0,\tau]$. The distribution of $H(\cdot)$ in $D_0[0,\tau]$ and the $k$-dimensional joint distribution of $(H(s_1),\ldots,H(s_k))$ under the probability $P_G$ will respectively be denoted by $\mathcal{L}(H;G)$ and $\mathcal{L}(H;G,s_1,\ldots,s_k)$. Let $\nu$ be a measure on $[0,\infty)$ with respect to which the distribution $G$ has a density $g$. Let $\mathcal{F}(\nu)$ denote the set of all densities with respect to $\nu$. Let $\mathcal{C}(g,\varsigma)$ be the set of all sequences of densities $\{g_k \in \mathcal{F}(\nu)\}$ such that

$$\lim_{k\to\infty} \|k^{1/2}(g_k^{1/2} - g^{1/2}) - \varsigma\|_2 = 0, \qquad (6.1)$$

where $\varsigma \in L_2(\nu)$ and $\|\cdot\|_2$ is the $L_2(\nu)$ norm. The limit in (6.1) implies in a standard way that $\varsigma \perp g^{1/2}$. Let $\mathcal{C}(g) = \bigcup_{\varsigma \in L_2(\nu), \varsigma \perp g^{1/2}} \mathcal{C}(g,\varsigma)$.

Suppose $\{g_k\} \in \mathcal{C}(g)$ is an arbitrary sequence with corresponding c.d.f.s $\{G_k\}$. Following Beran [5], we say that a sequence of estimators $\tilde{G}_k$ is regular at $g$ if

$$\mathcal{L}(\sqrt{k}(\tilde{G}_k - G_k); G_k) \xrightarrow{w} \mathcal{L}(\tilde{U}; G) \qquad \text{in } D_0[0,\tau],$$

where $\mathcal{L}(\tilde{U}; G)$ depends only on $g$ and not on the choice of the sequence $\{g_k\} \in \mathcal{C}(g)$ which determines the sampling scheme. Theorem 3 below establishes superiority of the NPMLE over all regular estimators whose finite-dimensional limiting laws are continuous. The proof is similar to the proof of Theorem 3 of [24] and, therefore, is omitted.

THEOREM 3. *Let $p > 0$ and $\tilde{G}$ be a sequence of regular estimators with a limiting law $\mathcal{L}(\tilde{U}; G)$ whose finite-dimensional laws $\mathcal{L}(\tilde{U}; G, s_1, \ldots, s_k)$ are continuous in $G$ under the supnorm topology for $G$. Then there exists a stochastic process $H(\cdot)$ in $D_0[0,t]$, $t \in \mathcal{J}$, such that*

$$\mathcal{L}(\tilde{U}; G) = \mathcal{L}(H; G) * \mathcal{L}(U; G),$$

*where $U$ is as in Theorem 2 and "$*$" denotes the convolution.*

**7. The other sampling schemes.** This section has two purposes: the first is to indicate briefly how the asymptotics might be established under scheme (ii), so that the case of random $M$ and fixed $M = m$ may be contrasted. The second is to demonstrate explicitly how multiplicative censoring [scheme (iii)] may be regarded as a special case of scheme (i).

*Sampling scheme* (ii): *Random censoring* (*conditional on* $M$). Under sampling scheme (ii), the proportion of uncensored observations is fixed. Assuming that $\hat{p} = p$, that is, conditioning on the proportion of uncensored observations, $\mathcal{H}_k$ given by (3.8) vanishes, $\mathcal{G}_{k,1}$ given by (3.6) reduces to $\mathcal{G}_1$ given by (5.2), while $\mathcal{G}_{k,2}$ given by (3.7) remains unchanged. Also, the second term on the left-hand side of (A.5) vanishes when we condition on the proportion of uncensored observations. We therefore obtain the following master equation for scheme (ii):

$$\Upsilon_k(U_{m,n}) = W_{m,n},$$

where

$$\Upsilon_k = \mathcal{G}_1 + \mathcal{G}_{k,2},$$

and $W_{m,n}$ is given by (3.3). It then follows from the results in Sections 4 and 5 that, under sampling scheme (ii), $\widehat{G}$ is uniformly strongly consistent and

$$U_{m,n} \overset{w}{\to} U = \mp^{-1}(W) \quad \text{in } D_0[0,t] \; \forall t \in \mathcal{J},$$

where

$$W(t) = p^{1/2} B_1(G_*(s)) + (1-p)^{1/2} f(s) \int_{0 < y \le s} B_2(F_*(y)) \, d\frac{1}{f(y)}$$

and $\Upsilon^{-1}$ is the inverse of $\Upsilon = \mathcal{G}_1 + \mathcal{G}_2$, where $\mathcal{G}_1$ and $\mathcal{G}_2$ are, respectively, given by (5.2) and (5.3).

*Sampling scheme* (iii): *Multiplicative censoring.* Having assumed (2.6) as the residual censoring distribution and by conditioning on the censoring proportion, $\mathcal{H}_k$ vanishes, while $\mathcal{G}_{k,1}(u)(t)$ and $\mathcal{G}_{k,2}(u)(t)$, respectively, reduce to $p_k I$ and $(1 - p_k) A_{\hat{f}}$, where $I$ is the identity map and

$$A_{\hat{f}}(u)(t) = \hat{f}(t) \int_{0 < y \le t} y \left( \int_{y \le z} \frac{u(z)}{z^2} \, dz \right) d\left(\frac{1}{\hat{f}(y)}\right).$$

Putting the above reduced forms together, we obtain the following master equation for sampling scheme (iii):

$$\Psi_k(U_{m,n}) = W_{m,n},$$

where $\Psi_k = p_k I + (1-p_k) A_{\hat{f}}$. It then follows from the results in Sections 4 and 5 that, under the multiplicative censoring scheme [scheme (iii)], $\widehat{G}$ is uniformly strongly consistent and

$$U_{m,n} \xrightarrow{w} U = \ominus^{-1}(W) \qquad \text{in } D_0[0,t],$$

where

$$W(t) = p^{1/2} B_1(G(s)) + (1-p)^{1/2} f(s) \int_{0<y\leq s} B_2(F(y)) \, d\frac{1}{f(y)},$$

and $\Psi^{-1}$ is the inverse of $\Psi = pI + (1-p)A_f$, if $p = \lim_{k\to\infty} p_k > 0.59$.

**8. Concluding remarks.** We have proved that, for length-biased right censored prevalent cohort survival data with follow-up, the NPMLE of the length-biased (i.e., prevalent case) survival function is strongly uniformly consistent, converges weakly to a Gaussian process, and is asymptotically efficient. It can be shown [4] that the NPMLE of the unbiased (i.e., incident case) survival function inherits these properties. The approach taken here is based on that used by Vardi and Zhang [24], although their methods do not carry over to the current more general setting without substantial modification, owing to the random censoring of the residual lifetimes. An apparently essential condition imposed for establishing the asymptotic results in Sections 4, 5 and 6 is that $t \in \mathcal{J}$. This condition is not restrictive since, in practice, $\beta$, the mass of the residual censoring distribution at 0, would be very small. For instance, if $\beta = 0.01$, then a sufficient condition for $t \in \mathcal{J}$ is that $F_C(t) \leq 0.98$.

In view of the fact that the current and residual lifetimes are equally distributed under stationarity, $\beta$ represents the proportion of uncensored observations with missing onset time. This then means that the results presented here address three of the four possible cases, that is, censored/uncensored and with/without onset time, in the setting considered in this paper.

If we allow an arbitrary unspecified incidence process, then the model becomes nonidentifiable and nonparametric estimation must be conditional on the truncation times, an approach that is commonly used because of its robustness against departure from stationarity of the incidence process. Wang [25], however, points out that this approach is only justified as conditional maximum likelihood if all censoring times are known, even for those who fail before they are censored. When the intensity of the incidence process is known, one can mimic the proofs given here to establish asymptotic results. This, however, entails a new master equation and, therefore, new subsequent steps.

## APPENDIX

PROOF OF LEMMA 1 (The master equation). Let $\Delta = [\hat{p}U_{m,n}(t) - W_{m,n}(t)]/\sqrt{k}$, where $\hat{p} = \frac{m}{k}$. Thus

$$\Delta = \hat{p}(\widehat{G}(t) - G_m(t)) - \frac{\sqrt{n}}{k}\hat{f}(t) \int_0^t W_{Y,n}(y)\, d\left(\frac{1}{\hat{f}(y)}\right) + \hat{p}(G_*(t) - G(t)).$$

Now using the equation (from [24], page 1034)

$$\frac{\hat{f}(t)}{\sqrt{n}} \int_{0<y\leq t} W_{Y,n}(y)\, d\left(\frac{1}{\hat{f}(y)}\right) = \int_{0<x\leq t}\int_{0<y\leq x} \frac{d(F_n(y) - F_*(y))}{\int_{y\leq z} z^{-1}\, d\widehat{G}(z)} x^{-1}\, d\widehat{G}(x),$$

obtained via integration by parts and a change of order of integration, we have

$$\Delta = \hat{p}(\widehat{G}(t) - G_m(t)) - (1-\hat{p})\int_{0<x\leq t}\int_{0<y\leq x} \frac{d(F_n - F_*)(y)}{\int_{y\leq z} z^{-1}\, d\widehat{G}(z)} x^{-1}\, d\widehat{G}(x)$$
$$+ \hat{p}(G_*(t) - G(t)).$$

Using (3.2),

$$\Delta = \hat{p}\widehat{G}(t) + (1-\hat{p})\int_{0<x\leq t}\left[\int_{0<y\leq x}\frac{dF_*(y)}{\int_{y\leq z} z^{-1}\, d\widehat{G}(z)}\right] x^{-1}\, d\widehat{G}(x)$$
$$- \widehat{G}(t) + \hat{p}(G_*(t) - G(t))$$
$$= (1-\hat{p})\int_{0<x\leq t}\left[\int_{0<y\leq x}\frac{dF_*(y)}{\int_{y\leq z} z^{-1}\, d\widehat{G}(z)} - dy\right] x^{-1}\, d\widehat{G}(x)$$
$$+ \hat{p}(G_*(t) - G(t)).$$

Utilizing (2.4), we obtain

$$\int_{0<y\leq x}\left[\frac{dF_*(y)}{\int_{y\leq z} z^{-1}\, d\widehat{G}(z)} - dy\right] = \int_{0<y\leq x}\left[\frac{F_C(y)f(y)\,dy/(1-p)}{\int_{y\leq z} z^{-1}\, d\widehat{G}(z)} - dy\right]$$
$$= \frac{1}{1-p}\int_{0<y\leq x}\left[F_C(y)\frac{f(y)}{\hat{f}(y)}\, dy - (1-p)\, dy\right].$$

We also have

$$\frac{F_C(y)f(y) - (1-p)\hat{f}(y)}{\hat{f}(y)}$$
$$= \frac{F_C(y)\int_{y\leq z} z^{-1}\, dG(z) - (1-p)\int_{y\leq z} z^{-1}\, d\widehat{G}(z)}{\int_{y\leq z} z^{-1}\, d\widehat{G}(z)}$$
$$= -\frac{\int_{y\leq z} z^{-1}\, dU_{m,n}(z)}{\sqrt{k}\int_{y\leq z} z^{-1}\, d\widehat{G}(z)} - \frac{S_C(y)\int_{y\leq z} z^{-1}\, dG(z)}{\int_{y\leq z} z^{-1}\, d\widehat{G}(z)} + p.$$

Thus,

$$\Delta = -\frac{(1-\hat{p})/(1-p)}{\sqrt{k}}$$

$$\times \int_{0<x\leq t}\left[\int_{0<y\leq x}\frac{\int_{z\geq y}z^{-1}\,dU_{m,n}(z)}{\int_{z\geq y}z^{-1}\,d\widehat{G}(z)}\,dy\right]x^{-1}\,d\widehat{G}(x)$$

$$+\left\{pG_*(t)\right.$$

(A.1)
$$-\frac{1-\hat{p}}{1-p}$$

$$\times \int_{0<x\leq t}\left[\int_{0<y\leq x}\frac{S_C(y)\int_{z\geq y}z^{-1}\,dG(z)}{\int_{z\geq y}z^{-1}\,d\widehat{G}(z)}\,dy\right]x^{-1}\,d\widehat{G}(x)\right\}$$

$$+p\left\{\frac{1-\hat{p}}{1-p}\widehat{G}(t)-G(t)\right\}$$

$$= I + II + III.$$

We simplify the terms $I$, $II$, $III$ in (A.1). First, as in [24], page 1035,

$$I = -\frac{(1-\hat{p})/(1-p)}{\sqrt{k}}\hat{f}(t)\int_0^t y\int_{z\geq y}\frac{U_{m,n}(z)}{z^2}\,dz\,d\left(\frac{1}{\hat{f}(y)}\right).$$

Next, in $II$, after substituting for $G_*(t) = \int_0^t g_*(x)\,dx$, using (2.3) and replacing $dG$ in the inner integral by $d\widehat{G} - k^{-1/2}\,dU_{m,n}$, we have

$$II = \int_0^t\left[\int_0^x\frac{\hat{p}-p}{1-p}S_C(y)\,dy\right]x^{-1}\,dG(x)$$

$$-\frac{1}{\sqrt{k}}\frac{1-\hat{p}}{1-p}\int_0^t\left[\int_0^x S_C(y)\,dy\right]x^{-1}\,dU_{m,n}(x)$$

$$+\frac{1}{\sqrt{k}}\frac{1-\hat{p}}{1-p}\int_0^t\left[\int_0^x\frac{S_C(y)\int_{z\geq y}z^{-1}\,dU_{m,n}(z)}{\int_{z\geq y}z^{-1}\,d\widehat{G}(z)}\,dy\right]x^{-1}\,d\widehat{G}(x),$$

while

$$III = p\left\{\frac{1-\hat{p}}{1-p}\frac{1}{\sqrt{k}}U_{m,n}(t) - \frac{\hat{p}-p}{1-p}G(t)\right\}.$$

Now, combining the above simplified forms for $I$, $II$, $III$, we obtain

$$\hat{p}U_{m,n}(t) - W_{m,n}(t)$$

$$= \left\{-\frac{1-\hat{p}}{1-p}\hat{f}(t)\int_0^t y\int_{z\geq y}\frac{U_{m,n}(z)}{z^2}\,dz\,d\left(\frac{1}{\hat{f}(y)}\right)\right\}$$

$$+ \left\{ \frac{\sqrt{k}(\hat{p}-p)}{1-p} \int_0^t \left[ \int_0^x S_C(y)\,dy \right] x^{-1}\,dG(x) \right.$$

$$- \frac{1-\hat{p}}{1-p} \int_0^t \left[ \int_0^x S_C(y)\,dy \right] x^{-1}\,dU_{m,n}(x)$$

$$\left. + \frac{1-\hat{p}}{1-p} \int_0^t \left[ \int_0^x \frac{S_C(y) \int_{z\geq y} z^{-1}\,dU_{m,n}(z)}{\int_{z\geq y} z^{-1}\,d\widehat{G}(z)}\,dy \right] x^{-1}\,d\widehat{G}(x) \right\}$$

$$+ p \left\{ \frac{1-\hat{p}}{1-p} U_{m,n}(t) - \frac{\sqrt{k}(\hat{p}-p)}{1-p} G(t) \right\}.$$

Thus,

(A.2)
$$\frac{\hat{p}-p}{1-p} U_{m,n}(t) + \frac{1-\hat{p}}{1-p}\hat{f}(t) \int_0^t y \int_{z>y} \frac{U_{m,n}(z)}{z^2}\,dz\,d\frac{1}{\hat{f}(y)}$$

$$+ \frac{1-\hat{p}}{1-p} \int_0^t \left[ \int_0^x S_C(y)\,dy \right] x^{-1}\,dU_{m,n}(x)$$

$$- \frac{1-\hat{p}}{1-p} \int_0^t \left[ \int_0^x \frac{S_C(y) \int_{z\geq y} z^{-1}\,dU_{m,n}(z)}{\int_{z\geq y} z^{-1}\,d\widehat{G}(z)}\,dy \right] x^{-1}\,d\widehat{G}(x)$$

$$= W_{m,n}(t) + p\sqrt{\frac{p}{1-p}}(G_*(t) - G(t)) \frac{\sqrt{k}(\hat{p}-p)}{\sqrt{p(1-p)}}.$$

Using the equation

$$\frac{d}{dy}\left( y \int_{z\geq y} z^{-2} U_{m,n}(z)\,dz \right) = -\int_{z\geq y} U_{m,n}(z)\,dz^{-1} - y^{-1} U_{m,n}(y)$$

$$= \int_{z\geq y} z^{-1}\,dU_{m,n}(z),$$

the fourth term on the left-hand side of (A.2) can be simplified to

(A.3)
$$\int_0^t \left[ \int_0^x \frac{S_C(y) \int_{z\geq y} z^{-1}\,dU_{m,n}(z)}{\int_{z\geq y} z^{-1}\,d\widehat{G}(z)}\,dy \right] x^{-1}\,d\widehat{G}(x)$$

$$= \int_0^t \left[ \int_y^t x^{-1}\,d\widehat{G}(x) \right] \frac{S_C(y)(d/dy)(y \int_{z\geq y} z^{-2} U_{m,n}(z)\,dz)}{\hat{f}(y)}\,dy$$

$$= \int_0^t \left( 1 - \frac{\hat{f}(t)}{\hat{f}(y)} \right) S_C(y)\,d\left( y \int_{y<z} z^{-2} U_{m,n}(z)\,dz \right)$$

$$= \int_0^t y \left( \int_{y<z} \frac{U_{m,n}(z)}{z^2}\,dz \right) d\left[ \left( \frac{\hat{f}(t)}{\hat{f}(y)} - 1 \right) S_C(y) \right].$$

Substituting (A.3) into (A.2), we obtain

$$\frac{\hat{p} - p}{1 - p} U_{m,n}(t) + \frac{1 - \hat{p}}{1 - p} \int_0^t y \left( \int_{z \geq y} \frac{U_{m,n}(z)}{z^2} \, dz \right) d\left[ \left( \frac{\hat{f}(t)}{\hat{f}(y)} - 1 \right) (1 - S_C(y)) \right]$$

(A.4) $\quad + \frac{1 - \hat{p}}{1 - p} \int_0^t \left[ \int_0^x S_C(y) \, dy \right] x^{-1} \, dU_{m,n}(x)$

$$= W_{m,n}(t) + p \sqrt{\frac{p}{1 - p}} (G_*(t) - G(t)) \frac{\sqrt{k}(\hat{p} - p)}{\sqrt{p(1 - p)}}.$$

Using (2.3) and (2.4), one can simplify (A.4) further to the form

$$\left( \frac{\hat{p} - p}{1 - p} \right) U_{m,n}(t)$$

$$+ \left[ 1 - \left( \frac{\hat{p} - p}{1 - p} \right) \right]$$

(A.5) $\quad \times \left\{ p \int_{0 < x \leq t} \frac{g_*(x)}{g(x)} \, dU_{m,n}(x) \right.$

$$\left. + (1 - p) \int_{0 < y \leq t} y \left( \int_{y \leq z} \frac{U_{m,n}(z)}{z^2} \, dz \right) d\left[ \left( \frac{\hat{f}(t)}{\hat{f}(y)} - 1 \right) \frac{f_*(y)}{f(y)} \right] \right\}$$

$$= W_{m,n}(t) + p \sqrt{\frac{p}{1 - p}} (G_*(t) - G(t)) \frac{\sqrt{k}(\hat{p} - p)}{\sqrt{p(1 - p)}}. \qquad \square$$

**Acknowledgment.** The authors thank an anonymous referee for his/her exceptionally careful reading of our manuscript, which resulted in the correction of several errors and the removal of some glaring ambiguities.

## REFERENCES

[1] ANDERSEN, P. K., BORGAN, Ø., GILL, R. D. and KEIDING, N. (1993). *Statistical Models Based on Counting Processes.* Springer, New York. MR1198884
[2] ASGHARIAN, M., WOLFSON, D. B. and ZHANG, X. (2004). A simple criterion for the stationarity of the incidence rate from prevalent cohort studies. Technical Report 2004-01, Dept. Mathematics and Statistics, McGill Univ.
[3] ASGHARIAN, M., WOLFSON, D. B. and ZHANG, X. (2005). Checking stationarity of the incidence rate using prevalent cohort survival data. *Statistics in Medicine.* To appear.
[4] ASGHARIAN, M., M'LAN, C. E. and WOLFSON, D. B. (2002). Length-biased sampling with right censoring: An unconditional approach. *J. Amer. Statist. Assoc.* **97** 201–209. MR1947280
[5] BERAN, R. (1977). Estimating a distribution function. *Ann. Statist.* **5** 400–404. MR0445701
[6] BICKEL, P. J., KLAASSEN, C. A. J., RITOV, Y. and WELLNER, J. A. (1993). *Efficient and Adaptive Estimation for Semiparametric Models.* Johns Hopkins Univ. Press, Baltimore. MR1245941

[7] BRILLINGER, D. R. (1986). The natural variability of vital rates and associated statistics (with discussion). *Biometrics* **42** 693–734. [MR0872958](MR0872958)

[8] CSÖRGŐ, M. and RÉVÉSZ, P. (1981). *Strong Approximations in Probability and Statistics.* Academic Press, New York. [MR0666546](MR0666546)

[9] FELLER, W. (1971). *An Introduction to Probability Theory and Its Applications* **2**, 2nd ed. Wiley, New York.

[10] GILBERT, P. B., LELE, S. R. and VARDI, Y. (1999). Maximum likelihood estimation in semiparametric selection bias models with application to AIDS vaccine trials. *Biometrika* **86** 27–43. [MR1688069](MR1688069)

[11] GILL, R. D., VARDI, Y. and WELLNER, J. A. (1988). Large sample theory of empirical distributions in biased sampling models. *Ann. Statist.* **16** 1069–1112. [MR0959189](MR0959189)

[12] GROENEBOOM, P. and WELLNER, J. A. (1992). *Information Bounds and Nonparametric Maximum Likelihood Estimation.* Birkhäuser, Basel. [MR1180321](MR1180321)

[13] HUANG, Y. and WANG, M.-C. (1995). Estimating the occurrence rate for prevalent survival data in competing risks models. *J. Amer. Statist. Assoc.* **90** 1406–1415. [MR1379484](MR1379484)

[14] KALBFLEISCH, J. D. and PRENTICE, R. L. (2002). *The Statistical Analysis of Failure Time Data*, 2nd ed. Wiley, New York. [MR1924807](MR1924807)

[15] KEIDING, N. (1990). Statistical inference in the Lexis diagram. *Philos. Trans. Roy. Soc. London Ser. A* **332** 487–509. [MR1084720](MR1084720)

[16] KEIDING, N., KVIST, K., HARTVIG, H., TVEDE, M. and JUUL, S. (2002). Estimating time to pregnancy from current durations in a cross-sectional sample. *Biostatistics* **3** 565–578.

[17] LEXIS, W. (1875). *Einleitung in die Theorie der Bevölkerungsstatistik*. Trübner, Strassburg. Pages 5–7 translated in (1977). *Mathematical Demography* (D. Smith and N. Keyfitz, eds.) 39–41. Springer, Berlin.

[18] LUND, J. (2000). Sampling bias in population studies—How to use the Lexis diagram. *Scand. J. Statist.* **27** 589–604. [MR1804165](MR1804165)

[19] PARTHASARATHY, K. R. (1967). *Probability Measures on Metric Spaces.* Academic Press, New York. [MR0226684](MR0226684)

[20] VAN ES, B., KLAASSEN, C. A. J. and OUDSHOORN, K. (2000). Survival analysis under cross-sectional sampling: Length bias and multiplicative censoring. *J. Statist. Plann. Inference* **91** 295–312. [MR1814785](MR1814785)

[21] VARDI, Y. (1982). Nonparametric estimation in the presence of length bias. *Ann. Statist.* **10** 616–620. [MR0653536](MR0653536)

[22] VARDI, Y. (1985). Empirical distributions in selection bias models (with discussion). *Ann. Statist.* **13** 178–205. [MR0773161](MR0773161)

[23] VARDI, Y. (1989). Multiplicative censoring, renewal processes, deconvolution and decreasing density: Nonparametric estimation. *Biometrika* **76** 751–761. [MR1041420](MR1041420)

[24] VARDI, Y. and ZHANG, C.-H. (1992). Large sample study of empirical distributions in a random-multiplicative censoring model. *Ann. Statist.* **20** 1022–1039. [MR1165604](MR1165604)

[25] WANG, M.-C. (1991). Nonparametric estimation from cross-sectional survival data. *J. Amer. Statist. Assoc.* **86** 130–143. [MR1137104](MR1137104)

[26] WANG, M.-C., BROOKMEYER, R. and JEWELL, N. P. (1993). Statistical models for prevalent cohort data. *Biometrics* **49** 1–11. [MR1221402](MR1221402)

[27] WANG, M.-C., JEWELL, N. P. and TSAI, W.-Y. (1986). Asymptotic properties of the product limit estimate under random truncation. *Ann. Statist.* **14** 1597–1605. [MR0868322](MR0868322)

[28] WOLFSON, C., WOLFSON, D., ASGHARIAN, M., M'LAN, C. E., ØSTBYE, T., ROCKWOOD, K. and HOGAN, D., for the Clinical Progression of Dementia Study Group (2001). A reevaluation of the duration of survival after the onset of dementia. *New England J. Medicine* **344** 1111–1116.

DEPARTMENT OF MATHEMATICS AND STATISTICS
MCGILL UNIVERSITY
BURNSIDE HALL
805 SHERBROOKE STREET WEST
MONTREAL, QUEBEC
CANADA H3A 2K6
E-MAIL: masoud@math.mcgill.ca
         david@math.mcgill.ca